\documentclass{amsart}
\usepackage{amsfonts,amsmath,amscd,amssymb,amsthm,newlfont}

\numberwithin{equation}{section}

\begin{document}

\title{A note on the invariance in the nonabelian tensor product}

%    Information for first author
\author{Francesco G. Russo}
%    Address of record for the research reported here
%\address{}
%    Current address
\curraddr{Laboratorio di Dinamica Strutturale e Geotecnica (StreGa)\\
Universit\'a del Molise, via Duca degli Abruzzi, 86039, Termoli
(CB).} \email{francescog.russo@yahoo.com}
%    \thanks will become a 1st page footnote.
\thanks{}

\date{\today}

%\dedicatory{This paper is dedicated to our advisors.}


\begin{abstract}
In the nonabelian tensor product $G\otimes H$ of two groups $G$ and $H$ many properties pass from $G$ and $H$ to
$G\otimes H$. There is a wide literature for different properties involved in this passage. We look at weak
conditions for which such a passage may happen.
\end{abstract}

\maketitle

\section{Terminology and statement of the result}

Let $G$ and $H$ be two groups acting upon each other in a
$compatible$ $way$:
\begin{equation}
~^{^gh}g'=~^{g}(^{h}(^{^{g^{-1}}}h')), \ \ \ \ \
~^{^hg}h'=~^{h}(^{g}(^{^{h^{-1}}}h')),
\end{equation}
for $g,g' \in G$ and $h,h' \in H$, and acting upon themselves by
conjugation. The $nonabelian$ $tensor$ $product$ $G\otimes H$ of $G$
and $H$ is the group generated by the symbols $g\otimes h$ with
defining relations
\begin{equation}
gg'\otimes h=(~^gg'\otimes ~^gh)(g\otimes h), \ \ \ \ g\otimes
hh'=(g\otimes h)(~^hg\otimes ~^hh').
\end{equation}
When $G=H$ and all actions are by conjugations, $G\otimes G$ is called $nonabelian$ $tensor$ $square$ of $G$.
These notions were introduced in \cite{BJR, BL} and some significant contributions can be found in \cite{bkm,
bmm, E, i, M1, M2, N, S, V}. %and the main calculus  rules in \cite[Propositions 1,2,3]{BJR}.
%\begin{equation}
% ~^g(g^{-1}\otimes h)=(g\otimes h)^{-1}=~^h(g\otimes h^{-1});
%\end{equation}
%\begin{equation}
% ~^{gh}(g'\otimes h')(g\otimes h)=(g\otimes
%h)~^{hg}(g'\otimes h');
%\end{equation}
%\begin{equation}
% g'\otimes (~^ghh^{-1})= ~^{g'}(g\otimes h)(g\otimes h)^{-1};
%\end{equation}
%\begin{equation}
% (g~^hg^{-1})\otimes h'=(g\otimes h)~^{h'}(g\otimes
%h)^{-1};
%\end{equation}
%\begin{equation}
% ~^{(g\otimes h)}(g'\otimes h^{'})=~^{[g,h]}(g'\otimes
%h');
%\end{equation}
%\begin{equation}
% [g\otimes h,g'\otimes h']=
%(g~^hg^{-1})(~^{g'}h'{h'}^{-1}).
%\end{equation}

From the defining relations in $G \otimes H$,
\begin{equation}
\kappa : g\otimes h \in G\otimes H\mapsto \kappa(g\otimes h)=[g,h] \in [G,H]=\langle g^{-1}h^{-1}gh \ | \ g\in
G, h \in H\rangle
\end{equation}
is an epimorphism of groups. Still from \cite{BJR, BL}, if $G$ and $H$ act trivially upon each other, then
$G\otimes H$ is isomorphic to the usual tensor product $G^{ab}\otimes_{\mathbb{Z}}H^{ab}$. If they act
compatibly upon each other, then their actions induce an action of the free product $G*H$ on $G\otimes H$ given
by $^x(g\otimes h) = ^xg \otimes ^xh$, where $x\in G*H$.

The $exterior$ $product$ $G\wedge H$ is the group obtained with the additional relation $g\otimes h=1_{\otimes}$
on $G\otimes H$, that is, \begin{equation}G\wedge H = (G\otimes H)/D,\end{equation} where $D=\langle g\otimes g
: g \in G\cap H\rangle$. Now it is easy to check that
\begin{equation}
\kappa': g\wedge h \in G\wedge H \mapsto \kappa'(g\wedge h)= [g,h]\in [G,H]\end{equation} is a well--defined
epimorphism of groups. For convenience of the reader, we recall that there is a famous commutative diagram with
exact rows and central extensions as columns in \cite[(1)]{BJR}: It correlates  the second homology group
$H_2(G)$ of $G$ with the third homology group $H_3(G)$ of $G$, the Whitehead's quadratic functor $\Gamma$, the
Whitehead's function $\psi$ and $\ker \kappa =J_2(G)$ (see also \cite{BJR, BL, W}).

%\begin{equation}
%\[\begin{CD}
%@. @.  0 @. 0\\
%@. @. @VVV   @VVV\\
%H_3(G)@>>>\Gamma(G^{ab}) @>\psi>>J_2(G)@>>>H_2(G)@>>>0\\
%@| @| @VVV @VVV\\
%H_3(G)@>>>\Gamma(G^{ab}) @>\psi>>G\otimes G @>>>G\wedge G@>>>1 @. \hspace{1.3cm}(*)\\
%@. @. @V\kappa VV   @V\kappa'VV\\
%@. @. G' @= G'\\
%@. @. @VVV   @VVV\\
%@. @. 1 @. 1\\
%\end{CD}\]
%\end{equation}

Now we get to the purpose of the present paper. Given a class of groups $\mathfrak{X}$, many authors answered
the question:
\begin{equation}\label{q:1}
\mathrm{If} \ \ G, H \in \mathfrak{X}, \ \ \mathrm{then} \ \  G\otimes H \in \mathfrak{X}\end{equation} In case
$\mathfrak{X}=\mathfrak{F}$ is the class of all finite groups, see \cite{E}. In case $\mathfrak{X}=\mathfrak{N}$
is the class of all nilpotent groups, see \cite{bmm, V}. In case $\mathfrak{X}=\mathfrak{S}$ is the class of all
soluble groups, see \cite{N, V}. In case $\mathfrak{X}=\mathfrak{P}$ is the class of all polycyclic groups, see
\cite{M1}. In case $\mathfrak{X}=\mathbf{L}\mathfrak{F}$ is the class of all locally finite groups, see
\cite{M2}. In case $\mathfrak{X}=\check{\mathfrak{C}}$ (resp., $\mathfrak{X}=\mathfrak{S}_2$) is the class of
all Chernikov (resp., soluble minimax) groups, see \cite{Rus}. Some topological properties are also closed with
respect to forming the nonabelian tensor product, as observed in \cite{BJR, BL}.

We recall some notations from \cite{LR}.
\begin{itemize}
\item[--]$\mathfrak{X}=\mathbf{S} \mathfrak{X}$ means that $\mathfrak{X}$ is closed with respect to forming subgroups.
\item[--]$\mathfrak{X}=\mathbf{H} \mathfrak{X}$ means that $\mathfrak{X}$ is closed with respect to forming homomorphic images.
\item[--] $\mathfrak{X}=\mathbf{P}\mathfrak{X}$ means that $\mathfrak{X}$ is closed with respect
to forming extensions, i.e.: if $N\in \mathfrak{X}$ is  a normal
subgroup of $G$ and $G/N\in \mathfrak{X}$, then $G\in \mathfrak{X}$.
\item[--] $\mathfrak{X}=\mathbf{H_2}\mathfrak{X}$ means that $\mathfrak{X}$ is closed with respect to
forming the second homology group, i.e.: if $G\in \mathfrak{X}$,
then $H_2(G)\in \mathfrak{X}$.
\item[--] $\mathfrak{X}=\mathbf{H_3}\mathfrak{X}$ means that $\mathfrak{X}$ is closed with respect to
forming the third homology group, i.e.: if $G\in \mathfrak{X}$, then
$H_3(G)\in \mathfrak{X}$.
\item[--] $\mathfrak{X}=\mathbf{T}\mathfrak{X}$ means that $\mathfrak{X}$ is closed with respect to
forming (usual) abelian tensor products , i.e.: if $A, B\in
\mathfrak{X}$ are abelian, then $A \otimes_{\mathbb{Z}} B\in
\mathfrak{X}$.
\end{itemize}

Our main contribution is below.

\medskip

\textbf{Main Theorem.} \textit{ Let $G$ and $H$ be two groups,  acting compatibly upon each other and
$\mathfrak{X}=\mathbf{S}\mathfrak{X}=\mathbf{H}\mathfrak{X}=\mathbf{P}\mathfrak{X}=\mathbf{H_2}\mathfrak{X}=\mathbf{H_3}\mathfrak{X}=\mathbf{T}\mathfrak{X}$.
If $G, H, \Gamma((G\cap H)^{ab}) \in \mathfrak{X}$, then $G\otimes H\in \mathfrak{X}$.}

\medskip

In \cite{bmm, E, M1, M2, N, Rus, V}, the  quoted results  follow from Main Theorem, when we choose
$\mathfrak{X}$ among $\mathfrak{F}, \mathfrak{N}, \mathfrak{S}, \mathfrak{P}, \mathbf{L}\mathfrak{F},
\mathfrak{\check{C}}, \mathfrak{S}_2$.

\section{Proof and some consequences}

We illustrate that it is possible to adapt an argument in
\cite[Section 2]{M1}.

\begin{proof}[Proof of Main Theorem]

Let $P=G*H/IJ$ be the Pfeiffer product of $G$ and $H$,
 where $I$ and $J$ are the normal closures in $G*H$ of
$\langle {^h}ghg^{-1}h^{-1}: g\in G, h \in H\rangle $ and $\langle
{^g}hgh^{-1}g^{-1}: g\in G, h \in H\rangle$, respectively. See
\cite{M1, W}. Note that $P$ is a homomorphic image of $G\ltimes
H$, hence $P \in \mathfrak{X}$. Here we have used
$\mathfrak{X}=\mathbf{H}\mathfrak{X}$. Let $\mu: G\rightarrow P$
and $\nu: H\rightarrow P$ be inclusions. Denote
$\overline{G}=\mu(G)$ and $\overline{H}=\nu(H)$. Then
$\overline{G}$ and $\overline{H}$ are normal subgroups of $P$ and
$P=\overline{G}\ \overline{H}$. Of course,  $\ker \mu \leq Z(G)$
and $\ker \nu \leq Z(H)$. An argument as in \cite[Proposition
9]{BJR} shows that the following sequence is exact:
\begin{equation}(G\otimes \ker \nu) \times (\ker \mu\otimes
H)\stackrel{i} \longrightarrow G\otimes H \longrightarrow
\overline{G}\otimes\overline{H} \longrightarrow 1,
\end{equation}
where $i$ is the inclusion $(g\otimes h',g'\otimes h)\mapsto
(g\otimes h')(g'\otimes h)$. It is easy to see that $\textrm{Im}~i
\leq Z (G \otimes H)$. Since $^{h}g=^{\nu(g)}g$ and
$^{g}h=^{\mu(g)}h$, $\ker \mu$ and $\ker \nu$ act trivially on $H$
and $G$, respectively.

Therefore,
\begin{equation}
G\otimes \ker \nu \simeq G^{ab}\otimes_{\mathbb{Z}} \ker
\nu^{ab}=G^{ab}\otimes_{\mathbb{Z}} \ker \nu
\end{equation} and
\begin{equation}\ker \mu \otimes H \simeq \ker \mu^{ab}\otimes_{\mathbb{Z}}
H^{ab}=\ker \mu\otimes_{\mathbb{Z}} H^{ab}.\end{equation}

In particular, $G\otimes \ker \nu \in \mathfrak{X}$.  Here we have
used $\mathfrak{X}=\mathbf{T}\mathfrak{X}$. Analogously, $\ker \mu
\otimes H \in \mathfrak{X}$. It follows that $\textrm{Im}~i \in
\mathfrak{X}$ because it is a homomorphic image of $(G\otimes \ker
\nu) \times (\ker \mu\otimes H) \in \mathfrak{X}$. Still we have
used $\mathfrak{X}=\mathbf{H}\mathfrak{X}$.

Since $\overline{G} \otimes \overline{H}\simeq (G\otimes
H)/\textrm{Im}~i$, it is enough to prove that
$\overline{G}\otimes\overline{H}\in \mathfrak{X}$. Here we have
used $\mathfrak{X}=\mathbf{P}\mathfrak{X}$. We may work with
$\overline{G}$ instead of $G$ and with $\overline{H}$ instead of
$H$ in order to get our result. Then there is no loss of
generality in assuming that $G$ and $H$ are normal subgroups of
$P$, $P=GH$, and all actions are induced by conjugation in $P$.
Note that $(G\wedge H)/\ker \kappa'$ is isomorphic to $[G,H]\leq
G\cap H \leq G \in \mathfrak{X}$ and so $(G\wedge H)/\ker
\kappa'\in \mathfrak{X}$. Here we have used
$\mathfrak{X}=\mathbf{H}\mathfrak{X}=\mathbf{S}\mathfrak{X}$. If
we prove $\ker \kappa' \in \mathfrak{X}$, then  $G\wedge H \in
\mathfrak{X}$ by $\mathfrak{X}=\mathbf{P}\mathfrak{X}$. If we
prove also $D\in \mathfrak{X}$, then $G\otimes H \in
\mathfrak{X}$, still by $\mathfrak{X}=\mathbf{P}\mathfrak{X}$ and
we are done.

By \cite[Theorem 4.5]{BL}, we have an exact sequence:
\begin{equation}\longrightarrow  H_3(P/G) \oplus
H_3(P/H)\longrightarrow \ker \kappa' \longrightarrow H_2(P)
\longrightarrow.
\end{equation} Since $P, P/G, P/H \in \mathfrak{X}$, we have $H_2(P), H_3(P/G), H_3(P/H) \in \mathfrak{X}$.
Here we have used
$\mathfrak{X}=\mathbf{H}\mathfrak{X}=\mathbf{H_2}\mathfrak{X}=\mathbf{H_3}\mathfrak{X}$.
On the other hand, $\ker \kappa'$ is an extension of $H_3(G/M)
\oplus H_3(G/N) \in \mathfrak{X}$ by $H_2(G) \in \mathfrak{X}$.
Therefore, $\ker \kappa' \in \mathfrak{X}$, as claimed.   Here we
have used $\mathfrak{X}=\mathbf{P}\mathfrak{X}$.

Having in mind the famous diagram \cite[(1)]{BJR}, it is easy to check that there exists a well--defined
homomorphism of groups $\psi: \Gamma((G \cap H)^{ab}) \rightarrow (G\cap H) \otimes (G\cap H)$. See
\cite[p.181]{BJR} or \cite{BL}. Then $\textrm{Im}~\psi =D\in \mathfrak{X}$, as claimed. Here we have used
$\mathfrak{X}=\mathbf{H}\mathfrak{X}$ and $\Gamma((G \cap H)^{ab}) \in \mathfrak{X}$.

The result follows.

\end{proof}

Note that $\Gamma(G^{ab})$ plays a fundamental role in deciding if $G\otimes G\in \mathfrak{X}$. This was
already noted in \cite[Section 3]{bmm} for the class of all free nilpotent groups of finite rank. Then it is
clear that the following corollary extends many results in \cite{bmm, E, M1, M2, N, Rus, V} in case of the
nonabelian tensor square.

\medskip

\textbf{Corollary.}

\textit{Assume $G=H$ in Main Theorem. If $G, \Gamma(G^{ab})\in
\mathfrak{X}$, then $G\otimes G\in \mathfrak{X}$.}

\medskip

We end with two observations on the invariance with respect to the nonabelian tensor product.

\medskip

\textbf{Remark 1.} Sometimes it is enough that $[G,H]\in \mathfrak{X}$ in order to decide whether $G\otimes H
\in \mathfrak{X}$. In case of $\mathfrak{X}=\mathfrak{N}$, or $\mathfrak{X}=\mathfrak{S}$, this can be found in
\cite{BJR, N, V}.

\medskip

The second deals with the \textit{universal property} of the nonabelian tensor products.  %We recall it to
%convenience of the reader. The function $\phi:G\times H\longrightarrow L$, where $L$ is an assigned group, is
%called a {\it crossed pairing} if
%\begin{equation}
%\phi(gg',h)=\phi(~^gg',~^gh)\phi(g,h)\ \ \ , \ \ \ \phi(g,hh')=\phi(g,h)\phi(~^hg,~^hh').
%\end{equation}

%\textbf{Proposition (See \cite{BL}).} \textit{Let $G$ and $H$ be two groups,  acting compatibly upon each other.
%Then there exists a group $A$, unique up to isomorphism, and a crossed pairing $\phi: G \times H \longrightarrow
%A$ such that for every group $L$ and each crossed pairing $\beta: G \times H \longrightarrow L$ there exists a
%unique group homomorphism $\overline{\beta}: A \longrightarrow L$ such that the following diagram is
%commutative:}

%\setlength{\unitlength}{10mm}
%\begin{picture}(8,2)(-4,-0.5)
%\put(.8,1){$G \times H$} \put(2.2,1){\vector(1,0){2}}\put(4.4,1){$A$}\put(4.3,.8){\vector(-1,-1){1}}
%\put(2.9,-.5){$L$}\put(1.8,.8){\vector(1,-1){1}}\put(3,1.15){$\phi$}
%\put(1.8,.1){$\beta$}\put(4.1,.1){$\overline{\beta}$} \put(5.1,.1){i.e. \ $\overline{\beta}\phi = \beta$.}
%\end{picture}

\medskip

\textbf{Remark 2.} In a certain sense the \textit{universal property of the nonabelian tensor products} (see
\cite{BL}) justifies Main Theorem, because it shows that we need at least
$\mathfrak{X}=\mathbf{S}\mathfrak{X}=\mathbf{H}\mathfrak{X}=\mathbf{P}\mathfrak{X}$, if we hope to answer
(\ref{q:1}) positively.

\bibliographystyle{amsalpha}

\end{document}